\documentclass[12pt]{amsart}

\usepackage{amsmath}
\usepackage{amssymb}
\usepackage{amsfonts}
\usepackage[latin1]{inputenc}

\numberwithin{equation}{section}
\newtheorem{theorem}{Theorem}[section]
\newtheorem{lemma}[theorem]{Lemma}
\newtheorem{definition}[theorem]{Definition}

\newtheorem{corollary}[theorem]{Corollary}
\newtheorem{fact}[theorem]{Fact}

\newtheorem*{maintheorem}{Main Theorem}
\newcommand{\nc}{\newcommand}
\nc{\cf}{{\rm cf}}\nc{\prob}{{\rm Prob}}
\nc{\RO}{{\rm RO}}
\nc{\TV}{{\rm TV}}
\nc{\potom}{\ensuremath{{\cal P}(\omega)}}
\nc{\potinf}{\ensuremath{[\omega]^\omega}}
\nc{\pfin}{\ensuremath{{\cal P}(\omega)/{\rm fin}}}
\nc{\potfin}{\ensuremath{[\omega]^{<\omega}}}
\nc{\proofend}{\makebox{} \hfill $\square$ \\}
\nc{\proofendof}[1]{\hspace*{\fill} $\square_{\rm #1}$ \\}
\newenvironment{myrules}
{\begin{list}{}
{
 \setlength{\leftmargin}{0.5in}
 \setlength{\labelwidth}{1cm}
 \setlength{\labelsep}{0.2in}
 \setlength{\parsep}{0.5ex plus 0.2ex minus 0.1 ex}
 \setlength{\itemsep}{0.3ex plus 0.2 ex minus 0ex}
}}{\end{list}}
\newcounter{subalph}
\newenvironment{subalph}{\begin{list}{(\alph{subalph})}%
 {\usecounter{subalph}
  \setlength{\leftmargin}{2em} \setlength{\labelwidth}{1.25em}
  \setlength{\labelsep}{0.75em} \setlength{\topsep}{0.5ex plus0.2ex minus0.2ex}
  \setlength{\parsep}{0pt}     \setlength{\itemsep}{\topsep}}}%
{\end{list}}
\newcount\skewfactor
\def\mathunderaccent#1#2 {\let\theaccent#1\skewfactor#2
\mathpalette\putaccentunder}
\def\putaccentunder#1#2{\oalign{$#1#2$\crcr\hidewidth
\vbox to.2ex{\hbox{$#1\skew\skewfactor\theaccent{}$}\vss}\hidewidth}}
\def\name{\mathunderaccent\tilde-3 }

\nc{\nname}{\name}
\nc{\englishfornameentry}{tilde} 
\nc{\such}{\, | \,}                   
\nc{\meager}{\ensuremath{{\cal M}}}
\nc{\lebesgue}{\ensuremath{{\cal N}}}
\nc{\gc}{\ensuremath{\frak c}}
\nc{\gd}{\ensuremath{\frak d}}
\nc{\gb}{\ensuremath{\frak b}}
\nc{\gt}{\ensuremath{\frak t}}\nc{\gh}{\ensuremath{\frak h}}
\nc{\gp}{\ensuremath{\frak p}}
\nc{\add}[1]{\mbox{\ensuremath{{\rm add}(#1)}}}
\nc{\cov}[1]{\mbox{\ensuremath{{\rm cov}(#1)}}}
\nc{\unif}[1]{\mbox{\ensuremath{{\rm unif}(#1)}}}
\nc{\cof}[1]{{\mbox{\ensuremath{\rm cof}(#1)}}}
\nc{\MA}{\mbox{\rm MA}}
\nc{\GCH}{\mbox{\rm GCH}}
\nc{\CH}{\mbox{\rm CH}}
\nc{\zfc}{\mbox{\rm ZFC}}
\newcommand{\cal}{\mathcal}
\newcommand{\nothing}[1]{}
\nc{\divs}{{c_0 \setminus \ell^1}}
\nc{\divser}{(\divs, \lessi^*)/\approy}
\nc{\bfin}{\RO(\pfin \setminus\{0\},\subseteq^*)}
\nc{\bdivser}{\RO(\divser)}
\nc{\inc}{{\rm INC}}
\nc{\com}{{\rm COM}}
\nc{\approy}{\makebox{}\!\!\approx}
\nc{\lessi}{\leqslant}
\nc{\gessi}{\geqslant}
\nc{\interior}[1]{{\rm int}(#1)}
\nc{\closure}[1]{{\rm cl}(#1)}
\nc{\PPP}{{P_{\delta,\omega_2}}}
\begin{document}

\title[On absolutely divergent series]
{On absolutely divergent series}

\author[Fuchino, Mildenberger, Shelah, Vojt\'a\v{s}]
{Saka\'e Fuchino\\ Heike Mildenberger\\ Saharon Shelah\\
Peter Vojt\'a\v{s}}

\thanks{The second author was partially supported by a Lise 
Meitner Fellowship of the State of North Rhine Westphalia} 

\thanks{The third author's research 
was partially supported by the ``Israel Science
Foundation'', administered by the Israel Academy of Science and Humanities.
This is  the third author's publication no.\ 593}

\thanks{The last author
was partially supported by the ``Alexander von Humboldt-Stiftung'', Bonn,
Germany and by grant 2/4034/97 of the Slovak Grant Agency}


\date{March 9, 1999}

\address{Saka\'e Fuchino, 
  Dept.\ of Computer Sciences,
  Kitami Institute of Technology,
  Koen-cho 165 Kitami,
  Hokkaido 090 Japan}

\email{fuchino@info.kitami-it.ac.jp, fuchino@math.fu-berlin.de}

\address{Heike Mildenberger,
Mathematisches Institut, Universit\"at Bonn,
Beringstr.~1,
53115 Bonn, Germany,
\and Mathematical Institute,
The Hebrew University of Jerusalem,
 Givat Ram,
 Jerusalem 91904, Israel}

\email{heike@math.uni-bonn.de, heike@math.huji.ac.il}

\address{
Saharon Shelah,
Mathematical Institute,
The Hebrew University of Jerusalem,
 Givat Ram,
Jerusalem 91904, Israel
}
 
\email{shelah@math.huji.ac.il}

\address{
Peter Vojt\'a\v{s},
Mathematical Institute,
Slovak Academy of Sciences,
Jesenn\'a~5,
04154 Ko\v{s}ice, Slovak Republic}
 
\email{vojtas@kosice.upjs.sk}

\begin{abstract}
We show that in the $\aleph_2$-stage countable support iteration of Mathias
forcing over a model of $\CH$ the 
complete Boolean algebra
generated by 
absolutely
 divergent series under eventual dominance is  not
 isomorphic to the completion of $P(\omega)/$fin.

This complements Vojt\'a\v{s}' result, that under $\cf(\gc) = \gp$ the two
algebras are isomorphic \cite{Vo93}.
\end{abstract}

\maketitle

\section{Introduction}

One of the traditional fields of real analysis is the study of 
asymptotic behaviour
of series and sequences; see e.g.\ the monographs of G.~H.~Hardy 
\cite{Hardy} and G.~M.~Fikh\-ten\-golz \cite{Fichtengolz}.
Among these topics
 is the classical problem of tests of absolute convergence and/or 
divergence of series of real numbers. Of specific importance is the 
comparison test, because many other tests, like Cauchy's (root) test, 
d'Alembert's (ratio)
test, and Raabe's test, are special instances of it.

\smallskip

We employ here a global point of view (implicit) of set theory, rather than 
looking at explicit series and tests (because these are only countably many
explicit ones, 
as our language is countable, and hence from a global point of view not very 
interesting). From this global --- set theoretic --- point of view  the study 
of comparison tests  is nothing else than the 
study of the ordering of eventual 
dominance on absolute values of the sequences, which describe the
entries that have to be summed up in a series, 
or on sequences with nonnegative entries, 
to which we restrict ourselves. A sequence $\bar{b}$ is eventually smaller 
than a sequence $\bar{a}$, denoted as $\bar{b} \lessi^* \bar{a}$,
if we have that $b_n\lessi a_n$ for all but finitely many $n$.

\smallskip

Note that the stronger information in the sense of convergence is 
carried by the eventually  greater sequences, in contrast to divergence 
where it is carried by the smaller ones. 
Hence we are interested in $\lessi^*$ on $\ell^1$ upwards, whereas 
on the set of divergent series
$c_0\setminus\ell^1$ the relation $\lessi^*$ is interesting downwards. 

\smallskip

There is a substantial difference between $(\ell^1, \gessi^*)$ and 
$(c_0\setminus\ell^1, \lessi^*)$, namely the first is directed 
and the second is not. 
For a directed ordering, questions about unbounded and dominating 
families are interesting. 
T.~Bartoszy\'nski \cite{Ba84} has shown that the minimum
size of an 
unbounded family of absolutely convergent series 
${\frak b} (\ell^1, \gessi^*)$ is equal to $\add{\lebesgue}$,
the additivity of the ideal of sets of Lebesgue measure zero. 
Dually, the minimal size of a dominating 
family ${\frak d} (\ell^1, \gessi^*)$
is $\cof{\lebesgue}$,
 the minimal size of a base of the ideal of sets of 
measure zero. This result says that in order to decide the absolute 
convergence of all series we need $\cof{\lebesgue}$ many series 
as parameters in the comparison test. This number is 
known to be consistently smaller than the size of 
the continuum $2^{\aleph_0}={\frak c}$.

\smallskip

On the opposite side, with divergence we need always continuum many 
divergent series as parameters for a comparison test in order to 
decide the divergence of all series. That is because there are continuum 
many incompatible divergent series below each divergent 
series. This observation 
together with the $\sigma$-closedness of $(c_0\setminus\ell^1, \lessi^*)$ 
raises 
the question what $(c_0\setminus\ell^1, \lessi^*)$ looks like from 
the Boolean theoretic point of view. In 
\cite{Vo93} P.~Vojt\'a\v s  has proved 
that the complete Boolean algebra generated by  
$(c_0\setminus\ell^1, \lessi^*)$ is isomorphic to the 
completion of the algebra
$\pfin$ of subsets of natural numbers equipped with eventual inclusion, 
assuming ${\frak p}=\cf({\frak c})$ (e.g.\ under $\CH$ or $\MA$). 
Moreover, T.~Bartoszy\'nski and 
M.~Scheepers \cite{BartoScheepers}
have shown that the ${\frak t}$-numbers of both orderings
are the same without additional hypotheses. 
This leads to the formulation of the problem
whether these two algebras are always isomorphic, 
in  all models of axiomatic set theory.

\smallskip

There is yet another striking phenomenon: F.~Hausdorff 
has shown (in \cite{Hausdorff}) that there is 
in $\zfc$ an $(\omega_1, \omega_1^*)$ gap in 
$(c_0, \lessi^*)$, such that the lower 
part of the gap consists of convergent series and the upper 
part consists of  divergent series. This is especially interesting when both 
$\add{\lebesgue}$ and ${\frak t}$ are greater than 
$\omega_1$. In this case we cannot 
approach the ``border between convergence and divergence'' from either
single side in $\omega_1$ steps, but we can do it in $\omega_1$ steps 
if we do it  simultaneously from both sides by a Hausdorff gap.

\smallskip 

To finish this introductory motivation, let us state that we can consider the 
classical study of asymptotic behaviour in the real 
analysis as a sort of study of 
forcing notions, because a better estimate and/or a
stronger result really corresponds to a
stronger forcing condition (in the case of non-directed orderings). 
Although it 
is historically a part of real analysis, it has gained new
interest, because of 
numerous applications in complexity theory in computer science. 

\smallskip

We consider the following complete Boolean algebras:

\smallskip

{\bf 1.} The {\bf algebra of regular
open sets in the partial order
$(\pfin \setminus \{0\}, \subseteq^*)$,}
called $\bfin$, where fin is the ideal of finite subsets of
$\omega$ and $\pfin$ is the set of all equivalence classes $a/{\rm fin}
= \{ b \in \potom \such b \triangle a \mbox{ is finite }\}$.
($a \triangle b = (a \setminus b) \cup (b \setminus a)$ 
is the symmetric difference of $a$ and $b$.)

We have that $a/{\rm fin} \subseteq ^* b/{\rm fin}$ iff
$a \subseteq^* b$, i.e.\ iff $a \setminus b$ is finite.
The element 0 is the class $\emptyset/{\rm fin}={\rm fin}$.

The partial order $P=(\pfin \setminus\{0\}, \subseteq^*)$ is separative,
i.e.\ 
$$\forall p,q \in P \left(p \not\lessi q \longrightarrow
\exists r \in P \; (r\lessi p \; \wedge \; r \perp q)\right),
$$
(where $r \perp q$ iff $\neg \exists s \:(s \lessi r \; \wedge
\; s \lessi q)$) 
or, in topological terms, for $p \neq q \in P$ we have that
$$
\interior{\closure{\{p' \such p' \lessi p \}}}
\neq \interior{\closure{\{q' \such q' \lessi q \}}},$$
where the interiors and closures are taken in the so-called cut topology
on $(P,\lessi)$, which is generated by the basic open sets
$\{ \{p' \such p' \lessi p \} \such p \in P \}$.
Hence the map $p \mapsto \interior{\closure{\{p' \such p' \lessi p \}}}$
is an embedding into the algebra of regular open subsets
of $P$, called $\RO(P)$.

\smallskip 

In general, for a partial order $(P,\lessi)$,
$A \subseteq P$ is called regular open iff
$$\interior{\closure{A}} =A.$$

As  shown in \cite[page~152]{Jech}, for any separative $(P,\lessi)$
 there is a unique complete Boolean algebra 
$\RO(P)$ into which --- leaving out
the Boolean algebra's zero element, of course ---
it can be densely embedded. 

\smallskip

{\bf 2.}
{\bf The algebra of regular open sets
$\bdivser$}, where
$\divs = \{ \bar{c} = \langle c_n \such n \in \omega \rangle \such
c_n \in {\mathbb R}_+ \wedge \lim c_n =0
\wedge \sum c_n = \infty \}$,
$\bar{d} \lessi^* \bar{c}$ iff for all 
but finitely many $n$ we have that $d_n \lessi c_n$.
This partial order $(\divs, \lessi^*)$ is not separative, see \cite{Vo94}. 
Hence
we take the separative quotient (see \cite[page 154]{Jech}):
We set $\bar{d} \approx \bar{c}$ iff
$\forall \bar{e} (\bar{e} \perp \bar{d} \leftrightarrow
\bar{e}\perp \bar{c})$.
Then we have that
$$(\bar{d}/\approy) \; (\lessi/\approy) 
\; (\bar{c}/\approy) \mbox{ iff }
\forall \bar{e} \lessi^* \bar{d} \; \bar{e} \not\perp \bar{c}.
$$
We write 
 $({\divs}, \lessi^*)/\approy$ for 
$({\divs}/\approy, \, \lessi^*\!/\approy)$, the
separative quotient, which is densely embedded into
$\bdivser$, the second object of our investigation. 

\smallskip

The purpose of this paper is to prove the following

\begin{maintheorem} \label{1.1}
In any extension got by the
 $\aleph_2$-stage countable support iteration of Mathias
forcing over a model of $\CH$, the complete Boolean algebra
generated by the separative quotient of absolutely
 divergent series under eventual dominance is  not
 isomorphic to the completion of $P(\omega)/$fin.
\end{maintheorem}

\smallskip

{\bf Notation and precaution:} 
We shall be using some partial orders as notions of
forcing as well. Then the stronger condition is the
{\em smaller} condition. Thus $\lessi$ in forcing will often coincide
with $\subseteq^*$ or  $\lessi^*$. For functions
$f,g \colon \omega \to {\mathbb R}$ we say $f \lessi^* g$ iff for
all but finitely many $n$, $f(n) \lessi g(n)$.
For subsets $A, B \subseteq \omega$ we write $A \subseteq^* B$ iff
$A \setminus B$ is finite.
The quantifier $\forall^\infty$ means ``for all but finitely many'',
and the quantifier $\exists^\infty$ means ``there are infinitely many''.
Names for elements in forcing extensions are written with
\englishfornameentry s under the object, like $\nname{x}$,
and names for elements of the ground model are written
with checks above the objects, like $\check{x}$.

\smallskip

Our notation follows Jech \cite{Jech} and Kunen \cite{Kunen}.
Recall that a subset $A $ of a partial order $(P,\lessi_P)$ is called
open iff it contains with any of its elements also all
stronger (i.e.\ $\lessi_P$ than the given element) conditions. 

If the ordering is clear, we shall often write only $P$ instead of 
$(P,\lessi_P)$ and $\lessi$ instead of $\lessi_P$.

\section{$\gh$-numbers}\label{S2}

The means to distinguish the two algebras are the $\gh$-numbers.
Therefore this section collects the facts we need about this
cardinal characteristic. Note that by a result of
Bartoszy\'nski and Scheepers \cite{BartoScheepers}
our two partial orders have the same $\gt$-numbers.
For information on $\gt$ and other cardinal characteristics we refer
the reader to \cite{vanDouwen}.

\begin{definition}\label{2.1}
\begin{subalph}
\item
A complete Boolean algebra $B$ is called $\kappa$-distributive iff
for every sequence of sets $\langle I_\alpha \such \alpha \in \kappa
\rangle$ and
every set $\{ u_{\alpha, i} \such i \in I_\alpha, \alpha \in \kappa \}$
of members of $B$ the equation
$$
\prod_{\alpha \in \kappa} \; \sum_{i \in I_\alpha} u_{\alpha, i}
= \sum_{f \in \prod_{\alpha \in \kappa} I_\alpha} 
\;\prod_{\alpha \in \kappa} u_{\alpha, f(\alpha)} $$
 holds.
\item For a partial order $(P, \lessi)$, $\gh(P,\lessi)$ 
is the minimal cardinal $\kappa$ such that
$\RO((P,\lessi)/\approy)$ is not $\kappa$-distributive.
If there is no such $\kappa$, let $\gh(P,\lessi)$ 
be undefined.
\item $\gh = \gh(\pfin \setminus \{0\}, 
\subseteq^\ast)$ is the well-known $\gh$-number
which was introduced by Balcar, Pelant and Simon in \cite{BPS}.
In fact, it could also be written 
$\gh = \gh(\potom \setminus {\rm fin}, \subseteq)$, 
since the separative quotient of 
$(\potom \setminus {\rm fin}, \subseteq)$ is $(\pfin \setminus \{0\}, 
\subseteq^\ast)$.
\end{subalph}
\end{definition}

The separative quotient of a separative order is (isomorphic to)
the order itself, 
and the set of regular open sets of a complete Boolean
 algebra (minus its zero) is (isomorphic to) the algebra itself.
Hence
\begin{equation}
\gh(P) = \gh(P/\approy) = \gh (\RO(P/\approy)).
\end{equation}

The following fact allows us to work with various equivalent
definitions of $\gh(P,\lessi)$.

\begin{fact}\label{2.2}
For any partial order $(P,\lessi)$  and cardinal $\kappa$ the 
following are equivalent:
\begin{myrules}
\item[(1)] $\RO((P,\lessi)/\approy)$
is $\kappa$-distributive.
\item[(2)] The intersection of $\kappa$ 
open dense subsets of $(P,\lessi)$
that are closed under $\approx$ 
is dense in  $(P,\lessi)$.
\item[(2')] The intersection of $\kappa$ 
open dense subsets of $(P,\lessi)/\approy$ 
is dense in  $(P,\lessi)/\approy$.
\item[(3)] Every family of $\kappa$ maximal antichains in 
$P$ has a refinement.
\item[(3')] Every family of $\kappa$ maximal antichains in 
$P/\approy$ has a refinement.
\item[(4)]
Forcing with $(P,\lessi)/\approy$ does not add a new 
function from $\kappa$ to ordinals.
\item[(5)] In the following game $G(P,\kappa)$ the player
\inc\ does not have  a winning strategy. 
The game $G(P,\kappa)$ is played in $\kappa$ 
rounds, and the two players \inc\ and \com\ 
choose $p_\alpha^\inc, p_\alpha^\com$ in the $\alpha$-th round such that
for all $\alpha < \beta < \kappa$,
$$
p_\alpha^\inc \gessi p_\alpha^\com 
\gessi p_\beta^\inc \gessi p_\beta^\com.
$$
In the end, player \inc\ wins iff the sequence of moves does not have 
a lower bound in $P$
or if at some round he/she does not have a legal move.
 Of course, \inc\ stands for ``incomplete'' and
\com\ stands for ``complete''.
\end{myrules}
\end{fact}

\proof 
The equivalence of (1) to (4) is well-known (even for not necessarily
separative partial orders!). We show:
a) that $\neg$(2) implies $\neg$(5) and 
b) $\neg$(5) implies $\neg$(3). This is also proved, for a different
game, where \com\ begins, and for a special Boolean algebra in
\cite{ShSp:494}. For $G(P,\omega)$, the equivalence of (2) and (5)
is also proved in \cite{Jech:distrBA}.

\smallskip

a) Suppose that we are given open dense sets $\langle
D_\alpha \such \alpha \in \kappa \rangle$ that are closed under
$\approx$ and  such that $A =
\bigcap\{D_\alpha \such \alpha \in \kappa \}$ is not dense. Define a 
winning strategy for \inc\ in $G(P,\kappa)$ as follows: 
For $\alpha \gessi 0$, \inc\ plays $p_\alpha^\inc \in D_\alpha$ 
such that $p_\alpha^\inc \lessi p_\beta^\com$ for all $\beta < \alpha$ and 
such that $A$ contains no element $\lessi p_\alpha^\inc$.
The first move is possible because $A$ is not dense.
This is clearly a winning strategy for \inc.

\smallskip

b) Let $\sigma$ be a winning strategy for $\inc$ in the game $G(P,\kappa)$.
We define maximal antichains $\langle A_\alpha \such \alpha \in \gamma \lessi
\kappa \rangle$ in $P$ such that if $\alpha < \beta < \gamma$
then $A_\beta$ is a refinement of $A_\alpha$ and if $p_\beta \in A_\beta$ 
and $p_\alpha \in A_\alpha$ is the unique member of $A_\alpha$ 
such that $ p_\alpha \gessi p_\beta$ then $\langle p_\alpha \such \alpha \in 
\beta \rangle$ are responses by $\sigma$ in an initial segment of a play,
i.e.,
\begin{align*}
\forall \alpha \lessi \beta \;\; \mbox{ for some } 
\langle p_\gamma^\com \such \gamma \in \alpha \rangle \;\; 
p_\alpha = p_\alpha^\inc =  \sigma (
\langle p_\gamma^\inc, p_\gamma^\com \such \gamma < \alpha \rangle).
\end{align*}

\smallskip

Suppose first that  $\langle A_\alpha \such \alpha \in \delta \rangle$ 
has been constructed. If the sequence does not have a refinement, then 
$\neg$(3) is proved. Otherwise suppose that there is some refinement
$B$ (which is of course, an antichain). Suppose that
$\delta = \delta' +1$. Then set
\begin{align*}
A'_{\delta} = \{\sigma(\langle p_\alpha^\inc, p_\alpha^\com \such 
\alpha \lessi \delta' \rangle) \such & \langle p_\alpha \such 
\alpha^\inc \lessi \delta' \rangle \mbox{ is decreasing through }\\
& \mbox{ all the $A_\alpha$, and $p_{\delta'}^\com \in B$,
and for $\alpha < \delta'$,}\\
&\mbox{$p_\alpha^\com$ is 
such that $p_\alpha^\inc \gessi
p_\alpha^\com \gessi p_{\alpha+1}^\inc$} \},
\end{align*}
and take $A_\delta \supseteq A'_\delta$ such that
$A_\delta$ is a maximal antichain 
If $\delta$ is a limit, then
\begin{align*}
A'_{\delta} = \{\sigma(\langle p_\alpha^\inc, p_\alpha^\com \such 
\alpha < \delta \rangle) \such & \langle p_\alpha \such 
\alpha^\inc \lessi \delta \rangle \mbox{ is decreasing through }\\
& \mbox{ all the $A_\alpha$, and for $\alpha < \delta$,
$p_\alpha^\com$ is}\\
&\mbox{such that $p_\alpha^\inc \gessi
p_\alpha^\com \gessi p_{\alpha+1}^\inc$}\},
\end{align*}
and again take for $A_\delta$ a maximal antichain containing $A'_\delta$.

If the construction did not stop before $\kappa$, 
then we would have found a $\lessi$-cofinal
part  $\langle
p_\alpha \such \alpha \in \kappa\rangle$ 
of a play  $\langle
p^\inc_\alpha, p^\com_\alpha \such \alpha \in \kappa\rangle$
according to $\sigma$ in which \inc\ loses, which would be a
contradiction.
\proofendof{\ref{2.2}}

\smallskip

\relax From Fact~\ref{2.2} we also get
\begin{corollary}\label{2.3}
The following are equivalent:
\begin{subalph}
\item \inc\ has a winning strategy in $G(P,\kappa)$.
\item \inc\ has a winning strategy in $G(P/\approy, \kappa )$.
\item \inc\ has a winning strategy in $G(\RO(P/\approy), \kappa )$.
\end{subalph}
\end{corollary}

\section{Distinguishing $\gh$-numbers; $\pfin$}\label{S3}

Complete Boolean algebras that are isomorphic have the same 
$\gh$-numbers. We use this obvious fact in order to derive our main theorem
from 
\begin{theorem}\label{3.1}
Let $G$ be generic for the
 $\aleph_2$-stage countable support iteration of Mathias forcing over a
model of $\CH$. Then  we have that in $V[G]$,
\begin{myrules}
\item[(a)] $\gh(\pfin, \subseteq^*) = \aleph_2$, and
\item[(b)] $\gh((\divs, \lessi^*)/\approy) = \aleph_1$.
\end{myrules}
\end{theorem}
\noindent{\em Beginning of proof.} 
We start with a ground model $V \models \CH$ and take an $\omega_2$-stage 
countable support iteration $P=\langle P_\alpha, \nname{Q}_\beta \such
\beta \in \omega_2, \alpha \lessi \omega_2 \rangle$ of Mathias forcing,
i.e.\ $\forall \alpha \in \omega_2$, 
$\Vdash_{P_\alpha} \mbox{``} \nname{Q_\alpha}$ is Mathias forcing''.

Remember that the conditions of Mathias forcing are pairs
$\langle u,A\rangle \in [\omega]^{<\omega} \times
\potinf $ such that $\max u < \min A$, ordered by $\langle v, B \rangle 
\lessi
\langle u,A \rangle$ iff
$u \subseteq v \subseteq u \cup A $ and $B \subseteq A$. Mathias forcing
will also (outside the iteration) be denoted by $Q_M$.

\smallskip

It is well-known (see \cite{ShSp:494}) that Mathias forcing can be decomposed
as $Q_M = Q'_M * \nname{Q''_M}$, where $Q'_M$ is $(\pfin \setminus\{0\},
\subseteq^*)$, which is
$\sigma$-closed  and adds as a generic a Ramsey ultrafilter $G'_M$,
and $\nname{Q''_M}$ denotes a name for
 Mathias forcing with conditions with second component
in $\nname{G'_M}$ (also know in the literature as 
${\mathbb M}_{\nname{G_M'}}$).
The ($Q'_M$-name for the) generic filter for $\nname{Q_M''}$ 
(which determines the Mathias real) will be denoted by $\nname{G_M''}$.
The map sending $\langle u,A \rangle$ to $\langle A, \langle u , A
\rangle \rangle$ is a dense embedding from
$Q_M$ into $Q'_M * \nname{Q''_M}$.

Since the first component is $\sigma$-closed and the second component is 
$\sigma$-centred (hence c.c.c.) the whole forcing is proper 
\cite{Properforcing}
and any iteration with countable support will not collapse $\aleph_1$.
Since for $\alpha < \omega_2$, $\Vdash_{P_\alpha} |\nname{Q}_\alpha
| \leq \omega_1$ and since the iteration length is $\lessi \omega_2$,
by \cite[III,4.1]{Properforcing}, $P$ has the $\aleph_2$-c.c.\ and hence
does not collapse any cardinals.

\smallskip

The next lemma is folklore. A proof of it with a slightly
more complicated argument can be found in  \cite{ShSp:494}.

\begin{lemma}\label{3.2}
In the model $V[G]$
from above we have that $\gh = \aleph_2$.
\end{lemma}

\proof Since in $V[G]$ we have that $2^\omega = \aleph_2$, we clearly have
$\gh \lessi \aleph_2$. We are now going to show that $\gh \gessi \aleph_2$.
We verify Fact~\ref{2.2}(2) for $\kappa= \aleph_1$.
In $V[G]$, let $\langle D_\nu \such \nu < \omega_1 \rangle$ be a family
of open dense sets of $\pfin \setminus\{0\}$.

By a L\"owenheim-Skolem 
argument, there is some $\omega_1$-club
(this is an unbounded set which is closed under suprema of strictly 
increasing $\omega_1$-sequences) $C \subseteq \omega_2$, $C \in V$, such that
for every $\alpha \in C$
$\forall \nu \in \omega_1$, $D_\nu \cap V[G_\alpha]$ is in $V[G_\alpha]$
and is open dense in $(\pfin)^{V[G_\alpha]} \setminus \{0\}$.
We want to prove that $\bigcap_{\nu \in \omega_1} D_\nu$ is not
empty below a given $B \in (\pfin)^{V[G]}
\setminus\{0\}$. By \cite{Properforcing}, there is
some $\delta < \aleph_2, \delta \in C$ such that $B \in V[G_\delta]$. 
By mapping $B$ bijectively, say via $f$, onto $\omega$ and changing 
the $D_\nu$ by mapping each of their members pointwise
with the same map $f$ we get $D_\nu'$, $\nu \in \omega_1$.
We claim the next Mathias real hits all the $D_\nu$ below
$B$. Now it is easy to see  that for $\nu \in \aleph_1$, that
\begin{equation*}
D_M(\nu) := \{ (u,A) \in Q_\delta \such 
A \in D'_\nu \cap V[G_\delta] \}
\end{equation*}
is dense in $Q_\delta$. 
So the Mathias real $r \in \potinf$ will be in all the $D_\nu'$. 
Now $f^{-1}{''} r $ is below $B$ and is in all the $D_\nu$.
\nothing{By genericity and a density argument for the
first components of the Mathias
forcings there is some $\alpha \in C$ such that $ A \in G(\alpha)'$,
where $G(\alpha)$ is the $Q_\alpha = 
\nname{Q}_\alpha[G_\alpha]$-generic filter 
determined by $G$ and $G_\alpha = G \cap P_\alpha$, 
and $G(\alpha)'$ is the first component of 
$G(\alpha)$ according to
the decomposition of the Mathias forcing explained above. 
Since $\alpha \in C$, $G(\alpha)'$ meets every $D_\nu$, $\nu \in \omega_1$.
But now $r_\alpha$, the $Q_\alpha''$-generic real determined by 
$(Q''_\alpha)_{Q_\alpha'}$, is below each member of $G(\alpha)'$, hence 
below $A$ and $\bigcap_{\nu \in \omega_1} D_\nu$. This proves that 
$\bigcap_{\nu \in \omega_1} D_\nu$ is dense.}
 \proofendof{\ref{3.2}}

\section{Distinguishing $\gh$-numbers; $\divs$}

In this section, we are going to prove 
$\gh((\divs,\lessi^*)/\approy)=\aleph_1$ in $V[G]$.
We work with the formulation
\ref{2.2}(2) and shall show something slightly stronger:

\begin{quote}
For any given $\bar{b} \in (\divs)^{V[G]}$, there are
$\langle D_\nu \such \nu \in \omega_1 \rangle
\in V[G]$ such that
$D_\nu$ is open and dense in $(\divs,\lessi^*)^{V[G]}$ and closed
under $\approx$
and such that their intersection is not dense below $\bar{b}$.
\end{quote}

\nothing{
We choose a family $\langle D_\nu \such \nu \in \omega_1 \rangle
\in V[G]$  such that 
$\langle D_\nu \such \nu \in \omega_1 \rangle$ is an enumeration of all 
 open dense subsets of $(\divs)^{V[G]}$ that are closed
under $\approx$ and are definable
with a real parameter $r \in (2^\omega)^{V}$.
Since the continuum hypothesis  
holds in $V$, such an enumeration exists.
}

Suppose that $\bar{b} \in (\divs)^{V[G]}$. There is some
$\delta < \omega_2$ such that $\bar{b} \in V[G_\delta]$.
We choose a family $\langle D_\nu \such \nu \in \omega_1 \rangle
\in V[G]$  such that 
$\langle D_\nu \such \nu \in \omega_1 \rangle$ is an enumeration of 
\begin{equation}\label{system}
\biggl\{
\Bigl\{ \bar{a} \in (\divs)^{V[G]} \Bigm| 
\sum_{\ell \in H} a_\ell <\infty \mbox{ or } \sum_{\ell \in 
\omega \setminus H} a_\ell <\infty \Bigr\} \biggm|
H \in (\potinf)^{V[G_\delta]} \biggr\}. 
\end{equation}
This is possible, because in $V[G_\delta]$ the continuum
has still cardinality $\aleph_1$.
All the sets in the set above are closed under $\approx$ and open and dense
in $(\divs)^{V[G]}$; the latter is shown as in Lemma~\ref{4.3}.
\smallskip

First let $\bar{m}(\bar{b}) =
\bar{m}= \langle m_i \such i \in \omega \rangle \in 
(\omega^\omega)^{V[G_\delta]}$ be a sequence of natural 
numbers such that for every $i \in \omega$,
\begin{eqnarray}
m_0 =0 \mbox{ and } m_{i+1} &>& 2^{m_i} \mbox{ and } \label{gross}\\ 
\frac{2^{m_{i}-2}}{36 \cdot (i+1)^2} &\gessi & 2^{(i+1)^2} \mbox{ and }\\
\forall \ell  \gessi m_{i+1} \;\; b_\ell &\lessi & 2^{-m_i}.\label{boundofb}
\end{eqnarray}

\nothing{Because of  \eqref{gross}, we have that
\begin{equation}
\label{divergent}
\sum_i \frac{m_{i+1} -m_i}{2^{m_{i}}} =  \infty.
\end{equation}
We set $b_\ell = \frac{1}{2^{m_i}}$ for $\ell \in [m_i,m_{i+1})$.
Then $\bar{b} = \langle b_\ell \such \ell \in \omega
\rangle \in \divs$ by \eqref{divergent}.
In the end we shall show that the intersection of the $D_\nu$
 is not dense
below $\bar{b}$.}

\smallskip

Now we begin an indirect proof. We assume 
\begin{equation}
\label{assumption}
\bigcap D_\nu \mbox{ is dense } (\lessi_\divs) \mbox{ below }
\bar{b}.
\end{equation}

\smallskip

The following chain of conclusions, including three lemmata,
 serves to derive a contradiction
from our assumption.
Following \cite{Baumgartner}, we factorise $P =
P_\delta * \PPP$.
We consider $V[G_\delta]$ as the ground model.
So there is a condition $p \in \PPP
\cap G$
and 
\begin{equation}
\label{dense}
p \Vdash_{\PPP} \mbox{``} \bigcap \nname{D_\nu}
\mbox{ is dense below } \check{\bar{b}} \mbox{''}.
\end{equation}

\smallskip

For technical reasons we have to ``discretize'' the partial
order $(\divs)^{V[G]}$ a bit. We set
\begin{multline*}
(\divs)^{\bar{m}}_{\rm discr} = \Bigl\{ \bar{e}=\langle e_\ell \such 
\ell \in \omega \rangle \in \divs
\Bigm| \forall i \in \omega \setminus\{0\}\; \forall \ell \in [m_i,m_{i+1})\\
          e_\ell \in \Bigl\{ \frac{j}{2^{m_{i+1}}} \Bigm| j =
 0,1,\dots , 2^{m_{i+1}-m_{i-1}} \Bigr\} \Bigr\}.
\end{multline*}

It is easy to see that $((\divs)^{\bar{m}}_{\rm discr})^{V[G]}$
(--- we interpret $\divs$ as a defining formula, which has to be evaluated
according to the model of set theory ---)
is dense in $(\divs)^{V[G]}$ below $\bar{b}$, the calculation that
$$
\sum_i \frac{m_{i+1} -m_i}{2^{m_{i+1}}} < \infty
$$
together with the formula~\eqref{boundofb}
helps to see it.

\smallskip

Because of $(\divs)^{\bar{m}}_{\rm discr}$'s
density below $\bar{b}$ and of \eqref{dense} we may assume that 
\begin{equation}
p\Vdash_{\PPP} \exists \bar{c} \lessi^*
\check{\bar{b}} \;\; \bar{c} \in (\divs)^{\bar{m}}_{\rm discr} 
\; \cap \;  \bigcap_{\nu \in
\omega_1} \nname{D_\nu}, 
\end{equation}
and we do so.

By the maximum principle, there is a name $\nname{\bar{c}}$
such that
\begin{equation}
\label{zufinale}
p\Vdash_{\PPP}   
\nname{\bar{c}} \in (\divs)^{\bar{m}}_{\rm discr} \cap \bigcap_{\nu \in
\omega_1} \nname{D_\nu} \; \wedge \; \nname{\bar{c}} 
\lessi^* \check{\bar{b}}.
\end{equation}
We set for $i \in \omega \setminus \{0\}$
$$
x_i = \left\{ s \Biggm| s \colon [m_i,m_{i+1}) \to 
\left\{ \frac{j}{2^{m_{i+1}}} \such j =  0,1,\dots, 2^{m_{i+1}-m_{i-1}} 
\right\} \right\}.
$$
Then we use

\begin{lemma}
\label{4.1}
(The Laver property for $\PPP$.) Suppose that
$\langle x_i \such i \in \omega \setminus \{0\}\rangle \in V[G_\delta]$ is a
family of finite sets and that
 $$
p \Vdash_{\PPP}
\forall i \in \omega \;\; \nname{\bar{c}} \!\restriction\! [m_i,m_{i+1})
 \in \check{x_i}.
$$
Then there are some $q \lessi_{\PPP} p$ and some
$\langle y_i \such i \in \omega \setminus \{0\}\rangle \in V[G_\delta]$
such that in $V[G_\delta]$
\begin{enumerate}
\item $\forall i \in \omega \setminus \{0\} \;\; 
|y_i| \lessi 2^{i^2}$, and
\item $\forall i \in \omega \setminus \{0\} \;\; y_i \subseteq x_i$, and
\item $q \Vdash_{\PPP}
\forall i \in \omega \setminus \{0\}\;\; \nname{\bar{c}} \!\restriction\! [m_i,m_{i+1})
 \in \check{y_i}.$
\end{enumerate}
\end{lemma}
\proof See Lemma 9.6. in \cite{Baumgartner}.

\medskip

Now we  apply Lemma~\ref{4.1} to our given $x_i$ and $\bar{c}$
and get
$\langle y_i \such i \in \omega \setminus \{0\}
\rangle \in V[G_\delta]$ as in the lemma.
We also fix some $q$ as in the lemma. Since there are
densely many such $q$ below $p$ and since $p \in G$ we may assume that
\begin{equation}\label{inG}
q \in G.
\end{equation}
For $i > 0$, we set 
$$
w_i = \Bigl\{ s \in y_i \Bigm| \sum_{\smash{\ell \in [m_i,m_{i+1})}}
s_\ell > \frac{1}{i^2} \Bigr\}.
$$

\smallskip

Since $\sum_{i \in \omega \setminus \{0\}} \frac{1}{i^2} < \infty$,
we have
that for any $\bar{e} \in (\divs)^{V[G]}$, 
\begin{multline}\label{bla}
\forall^\infty i \;\; \bar{e} \!\restriction\! [m_i,m_{i+1}) \in y_i 
\\
\longrightarrow 
\exists A \in \potinf \; (\forall i \in A \setminus \{0\}
 \;\; \bar{e} \!\restriction\! [m_i,m_{i+1}) \in w_i \; \wedge \;
\sum_{i \in A} \sum_{\ell \in [m_i,m_{i+1})} e_\ell = \infty).
\end{multline}

Note that by our choice of $\bar{m}$ we have for $i > 0$, 
\begin{equation}
|w_i| \lessi |y_i|
\lessi 2^{i^2} \lessi \frac{2^{m_{i-1}-2}}{36 \cdot i^2}.
\end{equation}

Before continuing in the main stream
of conclusions, we now record a useful lemma from probability theory.
The methods presented in \cite{AlonSpencer} led us to prove this lemma.

\begin{lemma}\label{4.2}
Assume that $\beta >0$, and 
\begin{myrules}
\item[(a)] $m < m' < m''$ are natural numbers, 
$m' > 2^m$ and $m'' > 2^{m'}$.
\item[(b)] $w \subseteq \left\{ \bar{d}
\such \bar{d} = \langle d_\ell \such m' \lessi \ell < m'' \rangle,
d_\ell \in \left\{ \frac{j}{2^{m''}} \such 0 \lessi j \lessi 2^{m''-m} 
\right\} \right\}$.
\item[(c)] If $\bar{d} \in w$, then 
$\sum\{ d_\ell \such \ell \in [m',m'')\} \gessi \frac{1}{\beta} $.
\item[(d)] $|w| \lessi \frac{2^{m-2}}{36 \beta}$.
\end{myrules}
Then we can find a partition $(u_0,u_1)$ of $[m',m'')$ such that
\begin{equation}\label{thirds}
\mbox{If $\bar{d} \in w$ and $h \in \{0,1\}$, }
\mbox{then } \frac{1}{3} \lessi \frac{\sum\{d_\ell \such \ell \in u_h \}}
{\sum \{ d_\ell \such \ell \in [m',m'')\}} \lessi \frac{2}{3}.
\end{equation}
 \end{lemma}

\proof We flip a fair coin for every $\ell \in [m',m'')$ to decide
whether $\ell$ is in
$u_0$ or in $u_1$ (so probabilities are $\frac{1}{2}$ and $\frac{1}{2}$).

\smallskip
We use $d = \sum_{\ell \in [m',m'')} d_\ell$ as an abbreviation.
Given $\bar{d} \in w$ and $h \in \{ 0,1\}$, we shall estimate
the probability

$$\prob\left( \frac{\sum\{ d_\ell \such \ell \in u_h \}}
{d}
< \frac{1}{3}\right).$$
The expected value of 
$$\frac{\sum\{ d_\ell \such \ell \in u_h \}}
{d}$$ 
is $\frac{1}{2}$.

\smallskip

TV denotes the truth value of an event $\varphi$:
$\TV(\varphi) = 1$ if $\varphi$ is true, and $\TV(\varphi)=0$ if
$\varphi$ is not true.
We compute the variance

\begin{align*}
{\rm Var} = & {\rm Exp}\left( \left(\frac{\sum\{ d_\ell \such \ell \in u_h \}}
{d} -
{\rm Exp}\left(\frac{\sum\{ d_\ell \such \ell \in u_h \}}
{d}\right) 
\right)^2\right)
\\
%
\intertext{ which equals, as the coins are thrown independently, }
=&  \frac{1}{d^2}
\cdot \sum_{\ell} 
\left( {\rm Exp}(d_{\ell}^2 \cdot {\rm TV}(\ell \in u_h))
- ({\rm Exp}(d_{\ell} \cdot {\rm TV}(\ell \in u_h)))^2 \right)\\
\lessi & \frac{1}{d^2}
\cdot \sum_{\ell} d_{\ell}^2 \cdot \frac{1}{2}.
\end{align*}

For the next argument, we allow, in contrast to our assumption (b)
of Lemma~\ref{4.2}, that 
the $d_\ell$ be reals such that
$$
0 \lessi d_\ell \lessi \frac{1}{2^{m}}.$$
We maximize 
$$\frac{1}{d^2}
\cdot \sum_{\ell} d_{\ell}^2 \cdot \frac{1}{2}
$$
 under the given requirements. The maximum of any variation
is attained if the
$d_\ell$, $\ell \in [m',m'')$ are most unevenly distributed, 
i.e.\ if some of 
them are $\frac{1}{2^{m}}$, one is possibly between 
0 and $\frac{1}{2^{m}}$ and the others are $0$.
In order to have them summed up to $d$, 
$v :=
\left\lfloor \frac{d}{\frac{1}{2^{m}}} \right\rfloor = \lfloor
2^{m} \cdot d \rfloor$ of them are $\frac{1}{2^{m}}$
(where $\lfloor x \rfloor$ 
denotes the  largest $n \in \omega$ such that $n \lessi x$).

\smallskip

Hence we get that
\begin{eqnarray*}
{\rm Var} &\lessi &  \frac{1}{d^2}
\cdot \sum_{\ell} d_{\ell}^2 \cdot \frac{1}{2}
\\
& \lessi & 
\frac{1}{ 2 \cdot {d^2}} \cdot 
\left( \left( \frac{1}{2^{m}}\right)^2 \cdot v + 
\left( (2^{m} \cdot d -v) \cdot  \frac{1}{2^{m}}\right)^2
\right) 
\\
& \lessi & \frac{1}{2\cdot d^2} \cdot 
\left( \frac{1}{2^{m}}\right)^2 \cdot d \cdot 2^{m}\\
&=& \frac{1}{2d} \cdot \frac{1}{2^{m}} \lessi 
\frac{\beta}{2^{m+1}} \mbox{ (see premise (c) of 
Lemma~\ref{4.2} for the last $\lessi$)}.
\end{eqnarray*}

We set
\begin{eqnarray*}
\alpha &=& \prob\left( \frac{\sum\{ d_\ell \such \ell \in u_h \}}
{d}
< \frac{1}{3}\right)\\
&=& \prob\left( \frac{\sum\{ d_\ell \such \ell \in u_h \}}
{d}
>\frac{2}{3}\right).
\end{eqnarray*}

So we get another estimate
\begin{eqnarray*}
\frac{\beta}{2^{m+1}} & \gessi & {\rm Var} = 
{\rm Exp} \left(\frac{\sum\{ d_\ell \such \ell \in u_h \}}
{d} - \frac{1}{2} \right)^2
\\
& \gessi & \alpha \cdot \left( \frac{-1}{6}\right)^2 +
\alpha \cdot \left( \frac{1}{6}\right)^2 = \frac{\alpha}{18}.
\end{eqnarray*}

Hence we have that
\begin{equation}
\label{bullet}
\alpha \lessi 18 \cdot \frac{\beta}{2^{m+1}}.
\end{equation}

The number of cases for a possible failure, which means 
 $\bar{d} \in w$ such that 
$$\frac{\sum\{ d_\ell \such \ell \in u_h \}}
{d}
\not\in \left[\frac{1}{3}, \frac{2}{3} \right], 
$$ 
is $|w|$, and the probability of any one failure is
$2\alpha$.

Hence we have at least one chance of success if
\begin{equation}\label{haha}
|w|\cdot 2\alpha <1, \end{equation}
because then
$$\prob(\mbox{no failure in $|w|$ cases}) \gessi 1- |w|\cdot 2 \alpha >0.$$
However, since by \eqref{bullet} we have that
$\alpha \lessi 18 \cdot \frac{\beta}{2^{m+1}}$,  and since
by our premises we have that $|w| \lessi \frac{2^{m-2}}{36 \beta}$,
our sufficient condition \eqref{haha} for success is fulfilled.
\proofendof{\ref{4.2}}

\smallskip

Now in $V[G_\delta]$ we apply Lemma~\ref{4.2} for every $i \in \omega$,
with $w = w_{i+1}$, $m = m_i$,  $m'=m_{i+1}$, $m'' =
m_{i+2}$, $\beta = (i+1)^2$,
 and we get for  $h =0,1$ for all $i \in \omega$ some 
 $u_{h,i+1} \subseteq [m_{i+1},m_{i+2})$ as in 
Lemma~\ref{4.2}.

\smallskip

With a real parameter in $V[G_\delta]$
(namely $\langle u_{0,i} \such i \in \omega \setminus \{0\} \rangle$)
we define the set
$$
J = \{ \bar{d} \in (\divs)^{V[G]} \such 
\exists h \in \{0,1\} \;\;\forall^\infty  i \in \omega \setminus \{0\}
\;\; \bar{d} \!\restriction\!
u_{h,i} \equiv 0 \}.
$$
$J$ is obviously open in $(\divs, \lessi^*)$. 

The closure of $J$ under $\approx$ is

\begin{equation}
\label{ugly}
J' = \{ \bar{d} \such \exists \bar{d'} \in J \;\;
 \forall \bar{e} \lessi^* \bar{d} \;\;
\; \bar{e} \not\perp \bar{d'} \}=
\{ \bar{d} \such \exists h \; \sum_{i \in \omega\setminus \{0\}
}\sum_{\ell \in u_{h,i}}
d_\ell < \infty\}.
\end{equation}

Note that we have 
\begin{equation}
\label{aboutJ}
(\bar{c}/\approy) \, \in  
\{ \bar{d}/\approy \such \bar{d} \in J \} \;\;\mbox{ iff }\;\;
\bar{c} \in J'.
\end{equation}
In the end, $J'$ will be the bad guy among the
$D_\nu$ from \eqref{zufinale}.

\smallskip

\begin{lemma}
\label{4.3}
$J$ is dense  in $(\divs)^{V[G]}$ under $\lessi^*$.
\end{lemma}

\proof Let $\bar{d}$ be an arbitrary 
element of  $(\divs)^{V[G]}$.
For $h \in \{0,1\}$ define $\bar{d}^h= \langle d_\ell^h \such \ell
\in \omega \rangle$ below $\bar{d}$ as follows
$$
d_\ell^h = \left\{\begin{array}{ll}
d_\ell, & \mbox{ if } \exists i \in \omega \setminus \{0\}\,( m_i \lessi  \ell 
< m_{i+1}, \mbox{ and } \ell \in
u_{h,i});\\
0, & \mbox{ else.}
\end{array}
\right.
$$
At least one of the $\bar{d}^h$ is divergent, because
$$
\sum d_\ell = \sum(d_\ell^0 +d_\ell^1).
$$
The divergent ones among the $\bar{d}^h$'s are in $J$.
\proofendof{\ref{4.3}}

Hence also $J'$ is dense.
So $J'$ is one of the $D_\nu$, namely with $H$ from~\eqref{system}
being $
\bigcup_{i \in \omega \setminus \{0\}} u_{0,i}$.
Now we can finally reach a contradiction by showing that
\begin{equation}
\label{finale}
q \not\Vdash_{\PPP} \nname{\bar{c}} \in \nname{J'} 
 \end{equation}
This will contradict \eqref{zufinale}.

\smallskip

In order to prove \eqref{finale}, we consider  
formula~\eqref{bla}, which yields
\begin{multline}
q \Vdash_{\PPP} 
\exists A \in \potinf \;\;(
\forall i \in A \setminus \{0\}\;\; \nname{\bar{c}}
\!\restriction\! [m_i,m_{i+1}) \in w_i
\; \wedge \; \sum_{i \in A} \sum_{\ell \in [m_i,m_{i+1})} \nname{c_\ell}
= \infty).
\end{multline}

Hence, by~\eqref{inG}, in $V[G]$ there is an infinite $A$ such that 
for $h =0,1$ we have by~\eqref{thirds} that
\begin{equation}\label{biggerthan}
V[G] \models 
\sum_{i \in A\setminus \{0\}} \sum_{\ell \in u_{h,i}} {c_\ell}
\gessi \frac{1}{3} \cdot 
\sum_{i \in A\setminus \{0\}} \sum_{\ell \in [m_i,m_{i+1})} {c_\ell} = \infty.
\end{equation}

Hence for either choice of $h \in \{0,1\}$ 
we have that $\bar{c}_h = \langle c^h_\ell \such \ell \in \omega \rangle$,
where
$$
c_\ell^h = \left\{\begin{array}{ll}
c_\ell, & \mbox{ if } \exists i \in A \setminus \{0\}\; (
 m_i \lessi  \ell < m_{i+1}, \mbox{ and } \ell \in
u_{h,i});\\
0, & \mbox{ else,}
\end{array}
\right.
$$
is divergent.

\smallskip

We shall show that $\bar{c} \not\in J'$ (though
$q \in G$), that is according to
the definition \eqref{ugly} of $J'$:
\begin{equation}
\forall \bar{d} \in J \;\; \exists \bar{c'} \lessi^* \bar{c}
\;\; \bar{c'} \perp \bar{d}.
\end{equation}

(Remark: Of course, we could have worked with
$\bar{c}/\approy$ and formulation \ref{2.2}(2$'$) all
the time and could have shown that there is  no $\bar{d} \in J$
that is $\approx \bar{c}$.
But we just did not like to handle equivalence classes all the time.)

\smallskip

Suppose we are given $\bar{d} \in J$. Then we have  that 
\begin{equation}
V[G] \models 
\exists h' \in \{0,1\} \;\; \forall^\infty i \in \omega \setminus \{0\} \;\;
{\bar{d}} \!\restriction\! u_{h',i} \equiv 0.
\end{equation}

We fix such a number $h'$.
But now we take $h =1-h'$! Then we have that $\bar{c}_{h}
\lessi^* \bar{c}$, and $\bar{c}_{h} $ is divergent, and
for every sequence $\bar{e}$ with
$(\bar{e} \lessi^* \bar{d} \; \wedge \; \bar{e} \lessi^*
 \bar{c}_h)$ we have that 
\begin{equation}
V[G] \models  \forall^\infty \ell \;\;  e_\ell = 0.
\end{equation}

Hence such an $\bar{e}$ cannot be a divergent series, and
we proved that $\bar{c}_h \perp \bar{d}$
and hence $\bar{c} \not\in J'$ (and, by \eqref{aboutJ},
$\bar{c} \not\approx \bar{d}$ for any
$\bar{d} \in J$).
This proves~\eqref{finale}.
So finally we derived a contradiction from~\eqref{assumption}.
\proofendof{{\rm Main Theorem}}

\smallskip

{\bf Acknowledgement:} The authors would like to thank 
Andreas Blass very much for carefully reading a preliminary version of
this paper, pointing out a gap, and making valuable suggestions.

\end{document}